\documentclass[letterpaper, 10 pt, conference]{ieeeconf}  

\IEEEoverridecommandlockouts                              

\usepackage{amsmath, amsfonts}
\usepackage{color}
\usepackage{amssymb}
\usepackage{mathtools}
\usepackage{csquotes}
\usepackage{quotes}
\usepackage{balance}

\usepackage{rotating}
\usepackage{adjustbox}

\usepackage{calc}

\usepackage{longtable}
\usepackage{multirow}

\usepackage{todonotes}

\usepackage{graphicx}

\usepackage{placeins}

\usepackage{tabularx}

\usepackage{lmodern,textcomp}

\usepackage{enumerate}

\usepackage{xspace}

\newcommand*\tageq{\refstepcounter{equation}\tag{\theequation}}

\usepackage{algorithm}
\usepackage{algpseudocode}

\newcommand{\inv}{^{-1}\xspace}

\newcommand{\blambda}{\boldsymbol{\lambda}\xspace}

\newcommand{\bxi}{\boldsymbol{\xi}\xspace}

\renewcommand{\t}{^\textsf{T}\xspace}

\newcommand{\away}[1]{}

\newcommand{\R}{\mathbb{R}\xspace}

\newcommand{\N}{\mathbb{N}\xspace}

\newcommand{\cL}{\mathcal{L}\xspace}

\newcommand{\cO}{\mathcal{O}\xspace}

\newcommand{\bx}{\textbf{x}\xspace}

\newcommand{\bg}{\textbf{g}\xspace}

\newcommand{\bd}{\textbf{d}\xspace}

\newcommand{\hbx}{\hat{\textbf{x}}\xspace}

\newcommand{\bz}{\textbf{z}\xspace}

\newcommand{\bei}[1]{{\textbf{e}}\xspace}

\newcommand{\bH}{\textbf{H}\xspace}
\newcommand{\bI}{\textbf{I}\xspace}

\newcommand{\bO}{\textbf{0}\xspace}

\newcommand{\tol}{{\textsf{tol}}\xspace}

\newcommand{\commentout}[1]{}

\usepackage{datetime,url}

\usepackage{cite}

\title{\LARGE \bf
A Direct Method for Solving 
Integral Penalty Transcriptions of Optimal Control Problems
}

\author{Martin P. Neuenhofen$^{1}$ and Eric C. Kerrigan$^{2}$
\thanks{$^{1}$Martin P. Neuenhofen is with the Deparment of Electrical \& Electronic Engineering, Imperial College London, SW7 2AZ London, UK, e-mail: m.neuenhofen19@imperial.ac.uk
        {\tt\small www.MartinNeuenhofen.de}}%
\thanks{$^{2}$Eric C. Kerrigan is with the Deparment of Electrical \& Electronic Engineering and Department of Aeronautics, Imperial College London, SW7 2AZ London, UK, e-mail: e.kerrigan@imperial.ac.uk
        {\tt\small www.imperial.ac.uk/people/e.kerrigan}}%
}

\begin{document}

\sloppy
\maketitle
\thispagestyle{empty}
\pagestyle{empty}

\begin{abstract}

We present a numerical method for the minimization of objectives that are augmented with large quadratic penalties of overdetermined inconsistent equality constraints. Such objectives arise from quadratic integral penalty methods for the direct transcription of equality constrained optimal control problems.
The Augmented Lagrangian Method (ALM) has a number of advantages over the Quadratic Penalty Method (QPM) for solving this class of problems. However, if the equality constraints of the discretization are inconsistent, then ALM might not converge to a point that minimizes the unconstrained bias of the objective and penalty term. Therefore, in this paper we explore a modification of ALM that fits our purpose.
Numerical experiments demonstrate that the modified ALM can minimize certain quadratic penalty-augmented functions faster than QPM, whereas the unmodified ALM converges to a minimizer of a significantly different problem.

\end{abstract}

\newcommand{\pval}{\omega\xspace} 	
\newcommand{\pmval}{\rho\xspace} 	
\newcommand{\wquad}{\alpha\xspace} 	
\newcommand{\tquad}{\tau\xspace} 	

\section{Motivation in the Optimal Control Context}

The method of choice for the numerical solution of optimal control problems is direct transcription. Typical direct transcriptions methods use orthogonal collocation~\cite{Betts2nd}. It is known that the latter can struggle with singular arc and high-index differential algebraic equalities (DAEs); the former arising in the example
\begin{equation*}
\begin{aligned}
&\min_{y,u} &\quad &J=\int_0^{\frac{\pi}{2}} \left(\,y(t)^2 + \cos(t)  u(t)\,\right)\,\mathrm{d}t,\\
&\text{s.t.} & y(0)&=0,\quad\dot{y}(t)=\frac{1}{2}  y(t)^2+u(t)\,,
\end{aligned}\tag{OCP}\label{eqn:ExampleOCP}
\end{equation*}
which has the analytic solution $y^\star(t)=\frac{\sin(t)}{\cos(t)-2}$ and $J^\star\approx -0.2569969625$\,.

Quadratic integral penalty methods \cite{Balakrishnan68,Hager90,PBF} are an alternative to collocation methods, where the squared path equality constraint residual is integrated and added as a penalty into the objective. In \cite{PBF} the authors present such a method with a  proof of convergence under mild assumptions, including convergence for singular arcs and high-index DAEs. This is verified in \cite{PBF} in comparison to collocation methods via numerical experiments.

Before proceeding, we guide the reader through the solution of \eqref{eqn:ExampleOCP} via the quadratic integral penalty method:
let $y,u$ be represented with continuous piecewise linear finite element functions $y_h,u_h$ on a uniform mesh of $N \in \N$ elements, $h:=\frac{\pi}{2 N}$; represented with $\bx := [y_h(h),\dots, y_h(Nh), u_h(0), \dots, u_h(Nh)]\t \in \R^n$, $n:={2 N+1}$. $y_h(0)=0$ is fixed and removed from $\bx$. We minimize
\begin{equation*}
\begin{aligned}
\min_{\bx \in \R^{2N+1}} \qquad &\Phi_{\pval}(\bx):=\int_0^{\frac{\pi}{2}} \left(\,y_h(t)^2 + \cos(t)  u_h(t)\,\right)\,\mathrm{d}t\\
+ \frac{1}{2  \pval}  &\int_0^{\frac{\pi}{2}} \left\|-\dot{y}_h(t)+\frac{1}{2} y_h(t)^2+u_h(t)\right\|_2^2\,\mathrm{d}t\,.
\end{aligned}\tageq\label{eqn:IntegralPenalty}
\end{equation*}
The integrals are evaluated with Gauss-Legendre quadrature of $q=8$ points per element. Writing $\tquad,\wquad$ for abscissae and weights, $m:=Nq$, and
\begin{subequations}
	\label{eqn:OCPdisc}
	\begin{align}
	f(\bx)&:= \sum_{j=1}^{Nq} \wquad_{j} \left( y_h(\tau_{j})^2 + \cos(\tau_{j})  u_h(\tau_j) \right)\label{eqn:OCPdiscF}\\
	c(\bx)&:= \begin{bmatrix} \vdots \\
	\sqrt{\wquad_{j}} \left( -\dot{y}_h(\tau_{j})^2 + \frac{1}{2}  y_h(\tau_j)^2 + u_h(\tau_j) \right)\\
	\vdots
	\end{bmatrix} \in \R^{m}\label{eqn:OCPdiscC}
	\end{align}
\end{subequations}
allows us to express \eqref{eqn:IntegralPenalty} as an \textit{unconstrained quadratic penalty program}:
\begin{equation*}
\begin{aligned}
\min_{\bx \in \R^{n}} \qquad &\Phi_{\pval}(\bx)= f(\bx) + \frac{1}{2  \pval}  \|c(\bx)\|_2^2\,,
\end{aligned}
\tag{UQPP}\label{eqn:UCQPP}
\end{equation*}
$\pval \in \R_+\setminus\lbrace 0\rbrace$ controls the quadratic penalty and should be chosen on the order of approximation of the finite element space \cite{Hager90,PBF}.

\newcommand{\bdual}{\blambda}
A sometimes related problem is the \textit{equality constrained program}:
\begin{equation*}
\begin{aligned}
\min_{\bx \in \R^{n}} \qquad &f(\bx) \qquad \text{s.t.}\quad c(\bx)=\bO \in \R^m\,,
\end{aligned}
\tag{ECP}\label{eqn:ECP}
\end{equation*}
with Lagrangian $\cL(\bx,\bdual) := f(\bx) - \bdual\t c(\bx)$ and Lagrange multiplier $\bdual \in \R^m$. The Karush-Kuhn Tucker optimality system of \eqref{eqn:UCQPP} and \eqref{eqn:ECP} is
\begin{equation*}
\begin{aligned}
\nabla_\bx \cL(\bx,\bdual)&=\bO\\
c(\bx) + \pval \bdual&=\bO
\end{aligned}
\tag{KKT}\label{eqn:KKT}
\end{equation*}
where $\pval=0$ for \eqref{eqn:ECP}. The KKT multiplier for \eqref{eqn:UCQPP} is a substitution trick such that $\nabla \cL$ matches~$\nabla \Phi_\omega$.

We recommend the use of the above penalty finite element method when numerically solving optimal control problems, because penalty methods have favourable convergence properties over collocation methods \cite{PBF}.

A difficulty that remains is with numerically solving \eqref{eqn:UCQPP}. This can be more challenging than solving a problem of the form \eqref{eqn:ECP}. Below we describe important details.

\subsection{Problems \eqref{eqn:UCQPP} and \eqref{eqn:ECP} have different solutions}
In our discretization \eqref{eqn:OCPdisc} it holds that $m=Nq\equiv8N$ and $n=2N+1\ll m$, hence \eqref{eqn:ECP} is \textit{overdetermined}.
Clearly, problem \eqref{eqn:UCQPP} cannot be overdetermined because it is unconstrained.

Considering problem \eqref{eqn:ECP} in the context of overdetermination poses the risk of \textit{inconsistency}. For instance, suppose that in \eqref{eqn:ExampleOCP} we had added the constraint $y(\pi/2)=-\frac{1}{2}$. Then $c(\bx)\neq\bO\ \forall \bx \in \R^n$ would hold,  i.e.\ \eqref{eqn:ECP} would be infeasible due to inconsistent constraints~$c$.
In contrast, problem \eqref{eqn:UCQPP} obviously possesses feasible points because it is unconstrained.

For our example, \eqref{eqn:ECP} is feasible,  i.e.\ $\exists \bx \in \R^n : c(\bx)=\bO$\,. Namely, the equality is achieved (only) by $\bx=\bO$, regardless of the discretization parameter~$h$. This is far away from $y^\star$, to which solutions of \eqref{eqn:UCQPP} converge as $h,\pval \searrow +0$. In conclusion, solutions to \eqref{eqn:ECP} can differ significantly from solutions to \eqref{eqn:UCQPP}.

The Modified Augmented Lagrangian Method (MALM), discussed in Section~\ref{sec:MALM}, converges to minimizers of \eqref{eqn:UCQPP} instead of \eqref{eqn:ECP}.

\subsection{Solutions of \eqref{eqn:UCQPP} depend on the value of $\pval$}
As experimentally verified in \cite{Hager90} and outlined in the analysis in \cite{PBF}, the discretization \eqref{eqn:IntegralPenalty}--\eqref{eqn:OCPdisc} converges when both $h,\pval \searrow +0$. That is, for fixed $h$, too large values of $\pval$ result in bad feasibility of the numerical optimal control solution, whereas too small values of $\omega$ result in feasible, yet far from optimal solutions.

\begin{figure}
	\centering
	\includegraphics[width=0.73\linewidth]{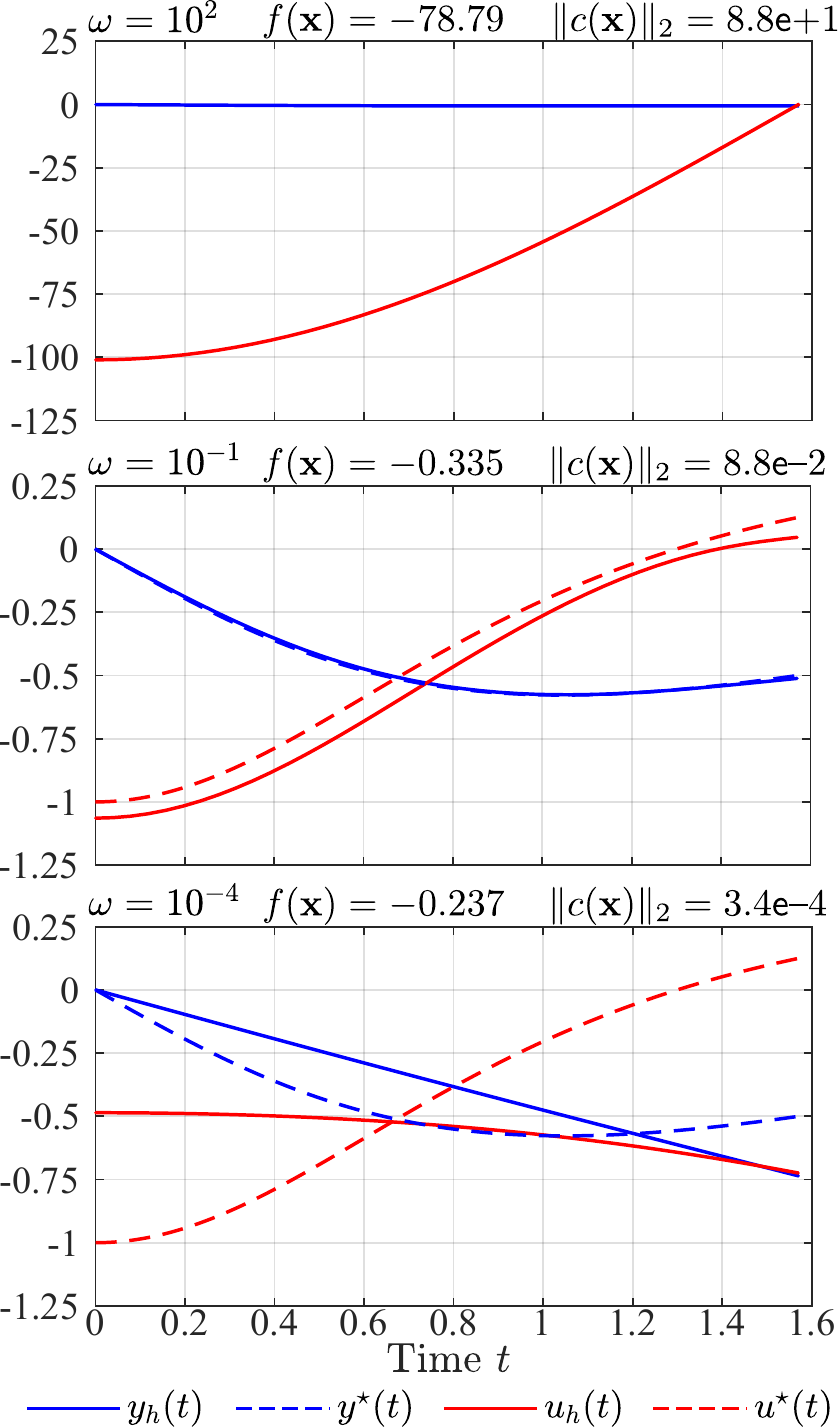}
	\caption{Numerical solution to \eqref{eqn:ExampleOCP} for $N=40$ and different values of $\pval$.}
	\label{fig:ocp_study}
\end{figure}

Fig.~\ref{fig:ocp_study} demonstrates this. Our discretization of \eqref{eqn:ExampleOCP} with $N=40$ is solved $\forall \pval \in \lbrace 10^{2},10^{-1},10^{-4}\rbrace$. The value~$\pval$ determines the bias between minimization of $f(\bx)$ and $\|c(\bx)\|_2 \equiv \|-\dot{y}_h+\frac{1}{2}y_h+u_h\|_{L^2(0,\pi/2)}$; note that the latter is accurate due to the order of quadrature. For $\pval=10^{-1}$ the solution achieves a good trade-off between feasibility and optimality on that coarse mesh.

In conclusion, the value of $\pval$ can have a significant influence on the solution of \eqref{eqn:UCQPP} and it is hence important that \eqref{eqn:UCQPP} be minimized for the specified value of $\pval$.

The value $\omega$ appears within the dual update formula of MALM. This is an important feature, so that MALM can converge to minimizers of \eqref{eqn:UCQPP} for the specific value of $\omega$.

\subsection{Direct minimization of \eqref{eqn:UCQPP} is numerically inefficient}
To the unprejudiced it appears natural to minimize the unconstrained objective $\Phi_{\pval}$ using a numerical method for unconstrained minimization. However, (unless $c$ is an affine function) this will result in very many iterations. This is so because the \textit{nonlinear} penalties result in \textit{curved} valleys through which iterative minimization methods make slow progress.

\newcommand{\abbTab}{Table}
\newcommand{\abbtab}{Table}
\newcommand{\abbAlg}{Algorithm}
\newcommand{\abbalg}{Algorithm}
\newcommand{\abbthm}{Theorem}
\newcommand{\abbThm}{Theorem}
\newcommand{\abbsec}{Section}
\newcommand{\abbSec}{Section}

To demonstrate this inefficiency, consider the instance
\begin{subequations}
	\label{eqn:Circle}
	\begin{align}
	f(\bx)&:= -x_1 -x_2\label{eqn:CircleF}\\
	c(\bx)&:= \begin{bmatrix}
		(x_1 +\varepsilon)^2 + x_2^2 - 2\\
		(x_1 -\varepsilon)^2 + x_2^2 - 2
	\end{bmatrix} \in \R^{m}\label{eqn:CircleC}
	\end{align}
\end{subequations}
with primal and dual initial guesses $\bx_0 :=\sqrt{2}[\cos(3\pi/8)\quad \sin(3\pi/8)]\t$ and $\bdual_0 :=0.4619[1\quad 1]\t$, for $\varepsilon=0$. We discuss later with \abbtab~\ref{tab:Exp_Circ_iter} that minimization of \eqref{eqn:UCQPP} of \eqref{eqn:Circle} with a simple trust-region method in \abbalg~\ref{algo:TRM} takes $334$ iterations when $\pval=10^{-6}$. This is inefficient when compared to MALM, which solves the same instance in only $16$ iterations.

The Augmented Lagrangiam Method (ALM) uses a local minimization method (typically a quasi-Newton variant) for the primal variables, and then updates the duals. Since ALM eventually solves \eqref{eqn:ECP}, the quasi-Newton system must solve \eqref{eqn:KKT} with $\omega=0$. In contrast, MALM converges to minimizers of \eqref{eqn:UCQPP}, i.e.\ solves \eqref{eqn:KKT} with $\omega>0$. This yields a dual regularization, which keeps the magnitude of $\|\bdual\|$ bounded and improves the convergence of the Newton iteration.

\subsection{Structure of the Paper}
Section II derives the proposed modified ALM (MALM) for nonlinear functions $f,c$.
Section III presents numerical experiments. This section also elaborates on the numerical difficulties of solving either \eqref{eqn:UCQPP} or \eqref{eqn:ECP} and suitable values of $\pval$ for a given instance $\bx_0,\bdual_0,f,c$. The numerical experiments compare the efficiency in terms of computational cost and iteration count for the Quadratic Penalty Method (QPM), ALM, and MALM.

\section{Derivation of the Modified Augmented Lagrangian Method}
\label{sec:MALM}

MALM is a solution method for \eqref{eqn:UCQPP}. MALM has been presented for minimizing convex 
quadratic penalty functions in \cite{SHARIFF2003257}, where $f$ must satisfy certain convexity properties and $c$ must be linear. Here, we derive MALM for nonlinear problems, in a stronger relation to its origins in ALM \cite{Hestenes1,Powell1}.

We derive MALM for \eqref{eqn:UCQPP} from ALM for \eqref{eqn:ECP}. To apply ALM, we need an auxiliary problem of form \eqref{eqn:ECP} instead. Our approach to achieving this works by temporarily using an auxiliary variable $\bxi \in \R^m$. This variable will be eliminated later in the augmented optimality system.

\subsection{Auxiliary Problem}
Consider the following equivalent problem to \eqref{eqn:UCQPP}:
\begin{subequations}
	\label{eqn:Subst}
	\begin{align}
	&\operatornamewithlimits{min}_{\hbx:=(\bx,\bxi) \in \R^{(n+m)}} 	&\quad 	&\hat{f}(\hbx):=f(\bx) + \frac{\pval}{2}  \|\bxi\|_2^2& \label{eqn:SubstObjective}	\\
	&\text{subject to} 						&		&\hat{c}(\hbx):=c(\bx) + \pval  \bxi =\bO\,. \label{eqn:SubstConstraints}
	\end{align}
\end{subequations}
Using Lagrange multipliers $\bdual \in \R^m$, the optimality conditions of \eqref{eqn:Subst} are
\begin{subequations}
	\label{eqn:KKT1_Subst}
	\begin{align}
	\begin{bmatrix}
	\nabla f(\bx)\\
	\pval\bxi
	\end{bmatrix} - \begin{bmatrix}
	\nabla c(\bx)\\
	\pval  \bI
	\end{bmatrix}  \bdual &=\bO\\
	c(\bx) + \pval  \bxi &= \bO\,. \label{eqn:5b}
	\end{align}
\end{subequations}

\subsection{Augmented Optimality System}
Since \eqref{eqn:Subst} is of form \eqref{eqn:ECP}, we can apply ALM \cite[Alg.~17.3]{NumOpt}. To this end, we augment \eqref{eqn:KKT1_Subst} with an auxiliary vector $\bz \in \R^m$ and a moderate penalty parameter $\pmval>0$:
\begin{subequations}
	\label{eqn:KKT1_Subst_ALM_unreduced}
	\begin{align}
	\begin{bmatrix}
	\nabla f(\bx)\\
	\pval\bxi
	\end{bmatrix} - \begin{bmatrix}
	\nabla c(\bx)\\
	\pval  \bI
	\end{bmatrix}  (\bdual+\bz) &=\bO\label{eqn:KKT1_Subst_ALM_unreduced1}\\
	c(\bx) + \pval  \bxi + \pmval  \bz &= \bO\,.\label{eqn:KKT1_Subst_ALM_unreduced2}
	\end{align}
\end{subequations}
The intuition for doing so in ALM is similar to what we did in \eqref{eqn:KKT} for \eqref{eqn:UCQPP}, where the Lagrange multipliers were used as a substitute to ensure that the gradient of~$\cL$ matches that of $\Phi_\omega$. Likewise, here $\bz$ works as a penalty substitute for \eqref{eqn:5b}.

We could use \eqref{eqn:KKT1_Subst_ALM_unreduced} directly in order to form an ALM iteration. That iteration would consist of alternately solving the optimality system \eqref{eqn:KKT1_Subst_ALM_unreduced} for $(\bx,\bxi,\bz)$ where $\bdual$ is fixed and updating $\bdual \leftarrow \bdual + \bz$, being equivalent to $\bdual \leftarrow \bdual - \frac{1}{\pmval}\left(c(\bx)+\pval\bxi\right)$.

\subsection{Elimination of the Auxiliary Vector}
However,  we propose to eliminate $\bxi = \bdual+\bz$ instead to obtain
\begin{subequations}
	\label{eqn:AlternatedKKT_MALM}
	\begin{align}
	\nabla f(\bx) - \nabla c(\bx)  (\bdual+\bz) &=\bO\\
	c(\bx) + \pval  \bdual + (\pval+\pmval)  \bz &= \bO \label{eqn:AlternatedKKT_MALM_primal}\,.
	\end{align}
\end{subequations}
As in ALM, we solve \eqref{eqn:AlternatedKKT_MALM} with an iteration of two alternating steps:
\begin{enumerate}
	\item Keep the value of $\bdual$ fixed, and solve \eqref{eqn:AlternatedKKT_MALM} for $(\bx,\bz)$.
	\item Update $\bdual$ as $ 	\bdual \leftarrow \bdual +\bz\,.	$
\end{enumerate}
Analogous to ALM, the first step can be realized by minimizing a suitable augmented Lagrangian function
for $\bx$ at fixed $\bdual$, whereas in the second step $\bz$ can be expressed in terms of $\bx$ from \eqref{eqn:AlternatedKKT_MALM_primal}. Using this, the method can be expressed compactly in \abbalg~\ref{algo:MALM}, where 
\begin{align}
\Psi_{k}(\bx) := \cL(\bx,\bdual_{k-1}) + \frac{1}{2 (\pval + \pmval)} \left\|c(\bx)+\pval  \bdual_{k-1}\right\|_2^2 \label{eqn:ALF}
\end{align}
is the augmented Lagrangian function, with $\cL(\bx,\bdual)=f(\bx)-\bdual\t \cdot c(\bx)$ as in \eqref{eqn:KKT}.

\begin{algorithm}[tb]
	\caption{Modified Augmented Lagrangian Method}
	\label{algo:MALM}
	\begin{algorithmic}[1]
		\Procedure{MALM}{$f,c,\pval,\bx_0,\bdual_0,\tol$}
		\State $\pmval \leftarrow \pmval_0$
		\For{$k=1,2,3,\dots,k_\text{max}$}
			\State Compute $\bx_k$ by TRM($\Psi_{k},\bx_{k-1},\tol$).
			$$\bx_k \leftarrow \operatornamewithlimits{argmin}_{\bx \in \R^n} \Psi_k(\bx) $$
			\State Update $\bdual_k \leftarrow \bdual_{k-1} - \frac{1}{\pval + \pmval} \left(c(\bx_k)+\pval \bdual_{k-1}\right)$
			\vspace{2mm}
			\If{$\|c(\bx_k)+\pval \bdual_k\|_\infty\leq\tol$}
				\State \Return $\bx_k,\bdual_k$
			\Else
				\State Decrease $\pmval\leftarrow c_{\pmval} \pmval$ to promote convergence.
			\EndIf
		\EndFor
		\EndProcedure
	\end{algorithmic}
\end{algorithm}

In our experiments we use $\tol=10^{-8}, c_\pmval=0.1, \pmval_0=0.1$. Care must be taken that the problem in line 4 is bounded below. To this end, practical methods use box constraints \cite[eq.~3.2.2]{Lancelot} or a trust-region \cite[eq.~3.2.4]{Lancelot}.

In order to minimize \eqref{eqn:ALF}, one can use any unconstrained minimization method. Here, for simplicity and reproducibility of the numerical experiments to follow, we use the simplified trust-region method in \abbalg~\ref{algo:TRM}. This method uses quasi-Newton directions $\bd$ in line~7, where the shift is determined directly by a decrease condition instead of implicitly by a trust-region radius. If we had used a line search method instead, then a shift would have still been necessary due to non-convexity, which would interfere with the line search; hence why we opted against it.

\begin{algorithm}[tb]
	\newcommand{\trnmobj}{\varphi}
	\caption{Simplified Trust-Region Method}
	\label{algo:TRM}
	\begin{algorithmic}[1]
		\Procedure{TRM}{$\trnmobj,\bx_0,\tol$}
		\State $j\leftarrow0$
		\While{$\|\nabla \trnmobj(\bx_j)\|_\infty>\tol$}
			\State $j\leftarrow j+1$,\qquad$\sigma\leftarrow 10^{-11}$
			\State $\bH\leftarrow\nabla^2 \trnmobj(\bx_{j-1})$,\qquad$\bg\leftarrow\nabla \trnmobj(\bx_{j-1})$
			\Repeat
				\State $\sigma\leftarrow10 \sigma$,\qquad$\bd\leftarrow-(\bH+\sigma\bI)\inv \bg$
			\Until{$\trnmobj(\bx_{j-1}+\bd)<\trnmobj(\bx_{j-1})$}
			\State $\bx_j \leftarrow \bx_{j-1}+\bd$
		\EndWhile
		\State \Return $\bx_j$
		\EndProcedure
	\end{algorithmic}
\end{algorithm}

As described in \cite[eq.~17.21]{NumOpt}, the quasi-Newton direction for the quadratic penalty function can be computed in a more numerically stable fashion from a $2 \times 2$ saddle-point linear equation system, by expressing the equations in terms of both $\bx$ and $\bz$. This has not been presented in the algorithms here for accessibility, but caused no issue in the numerical experiments with double precision and $\omega\geq 10^{-6}$.

\subsection{Discussion}
\subsubsection{True Generalization of ALM}\label{sec:truegeneralization}
MALM is a true generalization of ALM, because it differs merely by the parameter $\pval$. If $\pval=0$ then MALM is identical to ALM. Both methods then enjoy the same convergence properties and approach the same limit point.

\subsubsection{Benefit}
MALM solves the penalty function $\Phi_{\pval}$ in \eqref{eqn:UCQPP} by minimizing a sequence of penalty functions~$\Psi_k$. When does this make sense? By selecting $\pmval \gg \pval$. Thereby, the penalty functions~$\Psi_k$ have less steep valleys and hence can often be minimized more efficiently in comparison to one minimization of $\Phi_{\pval}$. The numerical experiments in the next section verify this claim.

\section{Numerical Experiments}
We present two numerical test problems. The first is instructional, the second is an optimal control problem.

\subsection{Circle Problem}
\subsubsection{Setting}
\begin{figure}
	\centering
	\includegraphics[width=0.58\linewidth]{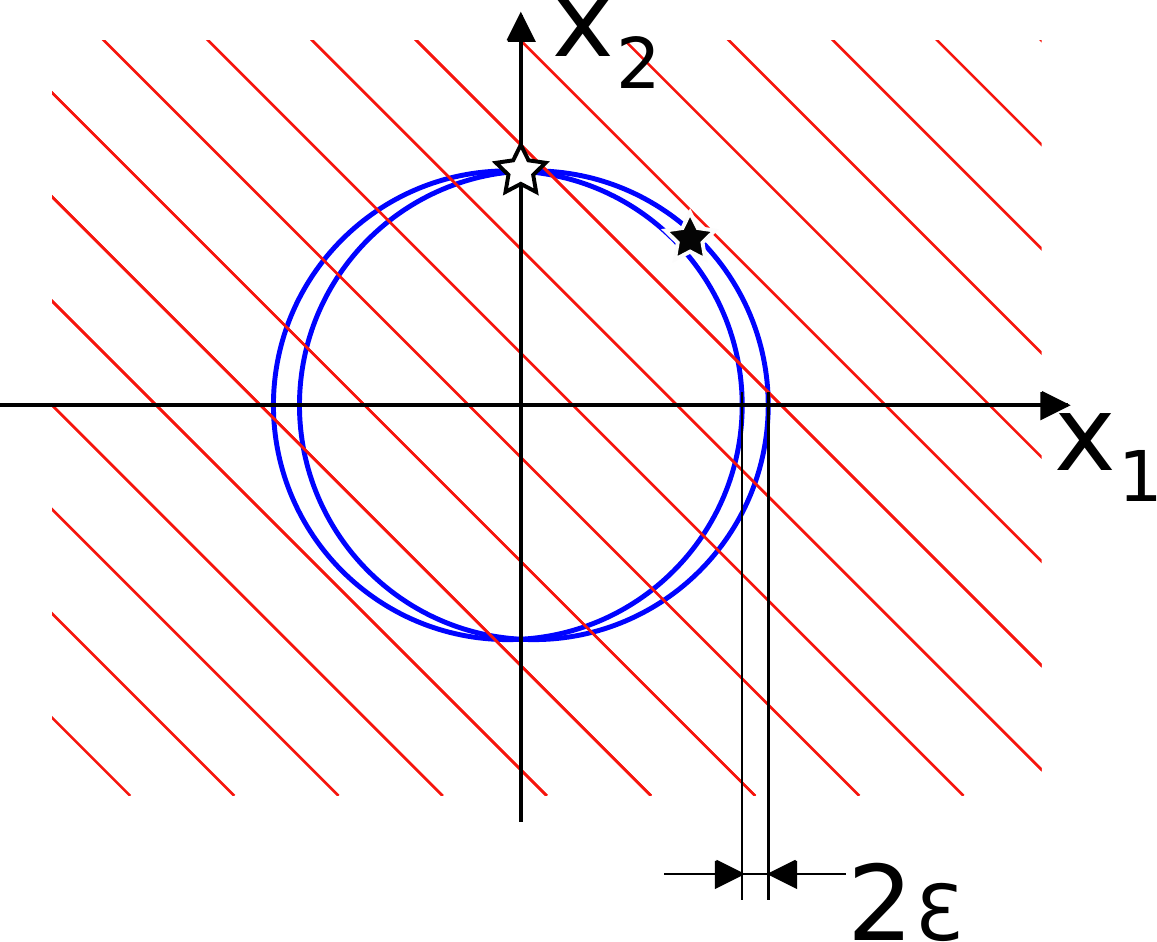}
	\caption{Level sets of components of $c$, and contours of $f$. The level sets of $c_1$ and $c_2$ are two circles, that almost fully overlap each other.}
	\label{fig:example1geometry}
\end{figure}

This problem considers the instance \eqref{eqn:Circle} for various values of $\varepsilon$. Fig.~\ref{fig:example1geometry} shows the geometry of the instance: Level sets of $f$ and $c$ are red and blue, respectively. The figure also shows two points
$$ 	\bx_A := [0\quad\sqrt{2}]\t\,,\qquad \bx_B := [1\quad1]\t 	$$
as the white and black star, respectively.

The instance can be interpreted in either of two ways:
\begin{enumerate}[$(A)$]
	\item Either we meant $c$ in a precise sense, meaning we wish to find a solution to $c(\bx)=\bO$ and, if non-unique, select the point that gives the smallest yield for $f$.
	\item Or we actually meant $c$ in a rough sense, meaning we wish to minimize $f$ subject to $\|\bx\|^2_2=2+\cO(\varepsilon)$.
\end{enumerate}

Both problems are reasonable in their own right: For example, $(A)$ makes sense when we have to solve a complex equation system and want to find a desirable solution. On the other hand, $(B)$ makes sense when our constraints suffer from errors, e.g.\ measurement errors or consistency errors, such as by discretization. For example, imagine a discretized optimal control problem, where~$c$ inherits consistency errors that have the size of $\varepsilon$.

Crucially, both solutions $\bx_A,\bx_B$ can be characterized sharply with a suitable problem statement. Obviously, $\bx_A$ is the solution of \eqref{eqn:ECP}. Less obvious, $\bx_B$ can be computed as the solution of \eqref{eqn:UCQPP} when choosing $\pval$ suitable w.r.t.\ $\varepsilon$. Here, a suitable choice is $\pval=\cO(\varepsilon)$. To see this, notice that \eqref{eqn:KKT} admits a well-scaled solution $\|\bx\|_2,\|\bdual\|_2 = \cO(1)$ and $\|c(\bx)\| = \cO(\varepsilon)$ when this selection for $\pval$ is made.

Lastly, we stress that for this instance the solution $\bx_A$ has an ill-conditioned KKT system with a dual solution $\|\bdual\| = \cO(1/\varepsilon^2)$, whereas $\bx_B$ is well-behaved, i.e, its KKT equations are well-conditioned.

\subsubsection{Computational Results}
We solve the instance with MALM and QPM, for various values of $\varepsilon,\pval$, including~0. We implement QPM by solving \eqref{eqn:UCQPP} via \abbalg~\ref{algo:TRM}. Recall that MALM=ALM for $\pval=0$ and that QPM is not applicable (n.a.) when $\pval=0$. 

Further, we investigate the limit points $\bx_\infty$ (which are identical for both tested methods throughout all tests) for each $\varepsilon,\pval$, by measuring the quantities
$$ 	e_A := \|\bx_\infty - \bx_A\|_2\,,\qquad e_B:=\|\bx_\infty-\bx_B\|_2\,.$$

\abbTab~\ref{tab:Exp_Circ_conv} shows the quantities $e_A,e_B$ for respective $\varepsilon,\pval$. Dividing the table into a lower left and an upper right triangle, we see that solutions in the lower triangle rather converge to $\bx_A$ while those on the diagonal and in the upper right converge to $\bx_B$.

\begin{table}[tb]
\caption{Solution of the Circle Problem w.r.t.\ $\varepsilon,\pval$. Smaller values mean closer convergence to either point. Table cells in the lower left converge to $\bx_A$, cells in the upper right to $\bx_B$.}
	\label{tab:Exp_Circ_conv}
	\centering
	\begin{tabular}{cl||c|c|c|c|c|}    \cline{3-7}\multicolumn{1}{l}{}&&\multicolumn{5}{c|}{$\varepsilon$}\\ \cline{2-7}  
		\multicolumn{1}{l|}{}& $\begin{matrix}e_A\\e_B\end{matrix}$& $1.0\text{e--}1$& $1.0\text{e--}2$& $1.0\text{e--}4$& $1.0\text{e--}6$& $0\text{e+}0$\\ \hline\hline
		\multicolumn{1}{|c|}{\multirow{5}{*}{$\pval$}} & $1.0\text{e--}1$&$\begin{matrix}6.5\text{e--}1\\4.6\text{e--}1\end{matrix}$&$\begin{matrix}1.1\text{e+}0\\1.0\text{e--}2\end{matrix}$&$\begin{matrix}1.1\text{e+}0\\8.8\text{e--}3\end{matrix}$&$\begin{matrix}1.1\text{e+}0\\8.8\text{e--}3\end{matrix}$&$\begin{matrix}1.1\text{e+}0\\8.8\text{e--}3\end{matrix}$\\ \cline{2-7}
		\multicolumn{1}{|c|}{\multirow{5}{*}{}} & $1.0\text{e--}2$&$\begin{matrix}1.2\text{e--}1\\9.7\text{e--}1\end{matrix}$&$\begin{matrix}1.0\text{e+}0\\5.6\text{e--}2\end{matrix}$&$\begin{matrix}1.1\text{e+}0\\8.8\text{e--}4\end{matrix}$&$\begin{matrix}1.1\text{e+}0\\8.8\text{e--}4\end{matrix}$&$\begin{matrix}1.1\text{e+}0\\8.8\text{e--}4\end{matrix}$\\ \cline{2-7}
		\multicolumn{1}{|c|}{\multirow{5}{*}{}} & $1.0\text{e--}4$&$\begin{matrix}3.7\text{e--}3\\1.1\text{e+}0\end{matrix}$&$\begin{matrix}1.1\text{e--}1\\9.8\text{e--}1\end{matrix}$&$\begin{matrix}1.1\text{e+}0\\5.7\text{e--}4\end{matrix}$&$\begin{matrix}1.1\text{e+}0\\8.8\text{e--}6\end{matrix}$&$\begin{matrix}1.1\text{e+}0\\8.8\text{e--}6\end{matrix}$\\ \cline{2-7}
		\multicolumn{1}{|c|}{\multirow{5}{*}{}} & $1.0\text{e--}6$&$\begin{matrix}3.5\text{e--}3\\1.1\text{e+}0\end{matrix}$&$\begin{matrix}1.2\text{e--}3\\1.1\text{e+}0\end{matrix}$&$\begin{matrix}1.0\text{e+}0\\5.6\text{e--}2\end{matrix}$&$\begin{matrix}1.1\text{e+}0\\5.7\text{e--}6\end{matrix}$&$\begin{matrix}1.1\text{e+}0\\8.8\text{e--}8\end{matrix}$\\ \cline{2-7}
		\multicolumn{1}{|c|}{\multirow{5}{*}{}} & $0\text{e+}0$&$\begin{matrix}3.5\text{e--}3\\1.1\text{e+}0\end{matrix}$&$\begin{matrix}3.5\text{e--}5\\1.1\text{e+}0\end{matrix}$&$\begin{matrix}3.5\text{e--}9\\1.1\text{e+}0\end{matrix}$&$\begin{matrix}3.5\text{e--}13\\1.1\text{e+}0\end{matrix}$&$\begin{matrix}1.1\text{e+}0\\0\text{e+}0\end{matrix}$\\ \hline
	\end{tabular}
\end{table}

\abbTab~\ref{tab:Exp_Circ_iter} shows the sum of the number of all inner iterations (i.e.\ iterations $j$ of \abbalg~\ref{algo:TRM}) of QPM and MALM for respective $\varepsilon,\pval$. We see a trend for each of the two methods: QPM converges in a few iterations when $\pval$ is moderate. However, when $\varepsilon,\pval$ both decrease, the iteration count blows up. The trend for MALM is different. MALM converges reliably for all $\varepsilon,\pval$ in the upper right triangle, including those where $\varepsilon,\pval$ are very small. The last row shows ALM. ALM converges quickly when $\varepsilon=0$, but its iteration count blows up for positive decreasing values of $\varepsilon$. In two instances ALM did not converge (n.c.) within $1000$ iterations.

\begin{table}[tb]
\caption{Number of iterations for MALM/ALM and QPM for the Circle Problem w.r.t.\ $\varepsilon,\pval$.  Fewer iterations mean better computational efficiency. MALM with $\omega=0$ (i.e.\ ALM) sometimes does not converge (n.c.) for this problem. QPM is not applicable (n.a.) when $\omega=0$.}
	\label{tab:Exp_Circ_iter}
	\centering
	\begin{tabular}{cl||c|c|c|c|c|}    \cline{3-7}\multicolumn{1}{l}{}&&\multicolumn{5}{c|}{$\varepsilon$}\\ \cline{2-7}  
		\multicolumn{1}{l|}{}&$\begin{matrix}\#_\text{MALM}\\\#_\text{QPM}\end{matrix}$& $1.0\text{e--}1$& $1.0\text{e--}2$& $1.0\text{e--}4$& $1.0\text{e--}6$& $0\text{e+}0$\\ \hline\hline
		\multicolumn{1}{|c|}{\multirow{5}{*}{$\pval$}} & $1.0\text{e--}1$& $\begin{matrix}\text{22}\\ \text{6}\end{matrix}$& $\begin{matrix}\text{19}\\ \text{9}\end{matrix}$& $\begin{matrix}\text{18}\\ \text{9}\end{matrix}$& $\begin{matrix}\text{16}\\ \text{9}\end{matrix}$& $\begin{matrix}\text{16}\\ \text{9}\end{matrix}$\\ \cline{2-7}
		\multicolumn{1}{|c|}{\multirow{5}{*}{}} & $1.0\text{e--}2$& $\begin{matrix}\text{25}\\ \text{9}\end{matrix}$& $\begin{matrix}\text{26}\\ \text{12}\end{matrix}$& $\begin{matrix}\text{19}\\ \text{13}\end{matrix}$& $\begin{matrix}\text{18}\\ \text{13}\end{matrix}$& $\begin{matrix}\text{16}\\ \text{13}\end{matrix}$\\ \cline{2-7}
		\multicolumn{1}{|c|}{\multirow{5}{*}{}} & $1.0\text{e--}4$& $\begin{matrix}\text{25}\\ \text{10}\end{matrix}$& $\begin{matrix}\text{83}\\ \text{30}\end{matrix}$& $\begin{matrix}\text{28}\\ \text{77}\end{matrix}$& $\begin{matrix}\text{21}\\ \text{77}\end{matrix}$& $\begin{matrix}\text{16}\\ \text{77}\end{matrix}$\\ \cline{2-7}
		\multicolumn{1}{|c|}{\multirow{5}{*}{}} & $1.0\text{e--}6$& $\begin{matrix}\text{26}\\ \text{10}\end{matrix}$& $\begin{matrix}\text{83}\\ \text{35}\end{matrix}$& $\begin{matrix}\text{97}\\ \text{311}\end{matrix}$& $\begin{matrix}\text{30}\\ \text{334}\end{matrix}$& $\begin{matrix}\text{16}\\ \text{334}\end{matrix}$\\ \cline{2-7}
		\multicolumn{1}{|c|}{\multirow{5}{*}{}} & $0\text{e+}0$& $\begin{matrix}\text{26}\\ \text{n.~a.}\end{matrix}$& $\begin{matrix}\text{83}\\ \text{n.~a.}\end{matrix}$& $\begin{matrix}\text{n.~c.}\\ \text{n.~a.}\end{matrix}$& $\begin{matrix}\text{n.~c.}\\ \text{n.~a.}\end{matrix}$& $\begin{matrix}\text{16}\\ \text{n.~a.}\end{matrix}$\\ \hline
	\end{tabular}
\end{table}

\subsubsection{Interpretation of the Results}
\abbTab~\ref{tab:Exp_Circ_conv} confirms that, depending on the parameters $\varepsilon,\pval$, we either solve for $\bx_A$ or $\bx_B$.
\abbTab~\ref{tab:Exp_Circ_iter} indicates that $\bx_A$ cannot be computed numerically efficiently, as expected due to the almost linearly dependent constraints.

In the optimal control context we are interested in solving problems in the sense of $(B)$, hence we now investigate the upper right triangle in \abbtab~\ref{tab:Exp_Circ_iter}. We see that both methods converge for all of these instances. Yet, when $\varepsilon,\pval$ are small, but strictly positive, then MALM outperforms QPM. This is relevant because there are problems from optimal control discretizations with inconsistencies like $\varepsilon$ from discretization errors. In these cases, we wish to drive $\pval,h\searrow+0$ to yield convergence of the discretization, which brings us into the lower right region of the table.

\subsection{Optimal Control Problem}
\subsubsection{Setting}
We solve the instance \eqref{eqn:OCPdisc} for various values of $N,\pval$ for $\bx_0=\bO,\bdual_0=\bO$. Recall that the instance represents the discretization \eqref{eqn:IntegralPenalty} with mesh size $h=\frac{\pi}{2 N}$, which only converges to the analytical solution when $h,\pval\searrow+0$ together, as was discussed along Fig.~\ref{fig:ocp_study}.

As for the former experiment, we solve the instance with MALM and QPM, for various values of $N,\pval$, including $\pval=0$. We also investigate the limit points $\bx_\infty$ (which are identical for both tested methods throughout all tests) for each $N,\pval$, by measuring the quantities
\newcommand{\symgap}{\delta J}
\begin{align*}
	\delta J 	&:= f(\bx)-J^\star\equiv J(y_h,u_h)-J^\star\\
	 	r 		&:=\|c(\bx)\|_2 \equiv \left\|-\dot{y}_h+\frac{1}{2} y_h^2+u_h\right\|_{L^2(0,\pi/2)}\,.
\end{align*}

Recall from the circle test problem that there are two interpretations $(A)$ and $(B)$ for the instance. As is clear from the context, we wish to compute a solution of kind $(B)$. However, if we are (deliberately) uncareful in choosing $\pval$ suitable w.r.t.\ $h$ then the iteration will converge to a solution of kind $(A)$.

\subsubsection{Computational Results}
\abbTab~\ref{tab:Exp_OCP_conv} shows the quantities $\symgap,r$ for respective $N,\pval$. Dividing the table into a lower left and an upper right triangle, we see that solutions in the lower left achieve good feasibility but at the sacrifice of optimality, whereas solutions in the upper right are not sufficiently feasible. Only solutions on the diagonal strike a good balance between minimizing $J$ and~$r$,  i.e.\ optimality and feasibility. Thus, when being limited by computation time to solve on a moderately sized mesh, then accordingly $\omega$ should not be chosen too small. Therefore, we must consider the table column by column. The table clearly shows that the best balance between feasibility and optimality is obtained in those table entries that live on the diagonal.

\begin{table}[tb]
\caption{Solution of the Optimal Control Problem w.r.t.\ $N,\pval$. For a given mesh size $N$, the value for $\omega$ is suitable when it yields a good balance between the orders of magnitude for both $\delta J$ (optimality gap) and $r$ (feasibility residual). Table cells on and around the diagonal typically yield a good balance.}
	\label{tab:Exp_OCP_conv}
	\centering
	\begin{tabular}{cc||c|c|c|c|c|}    \cline{3-7}\multicolumn{1}{l}{}&&\multicolumn{5}{c|}{$N$}\\ \cline{2-7}  
		\multicolumn{1}{l|}{}& $\begin{matrix}\symgap\\r\end{matrix}$& $16$& $64$& $256$& $1024$& $4096$\\ \hline\hline
		\multicolumn{1}{|c|}{\multirow{5}{*}{$\pval$}} & $1.0\text{e--}1$&$\begin{matrix}\text{-}7.8\text{e--}2\\9.2\text{e--}2\end{matrix}$&$\begin{matrix}\text{-}7.9\text{e--}2\\8.9\text{e--}2\end{matrix}$&$\begin{matrix}\text{-}7.9\text{e--}2\\8.9\text{e--}2\end{matrix}$&$\begin{matrix}\text{-}7.9\text{e--}2\\8.9\text{e--}2\end{matrix}$&$\begin{matrix}\text{-}7.9\text{e--}2\\8.9\text{e--}2\end{matrix}$\\ \cline{2-7}
		\multicolumn{1}{|c|}{\multirow{5}{*}{}} & $2.5\text{e--}2$&$\begin{matrix}\text{-}1.5\text{e--}2\\2.7\text{e--}2\end{matrix}$&$\begin{matrix}\text{-}2.0\text{e--}2\\2.3\text{e--}2\end{matrix}$&$\begin{matrix}\text{-}2.0\text{e--}2\\2.2\text{e--}2\end{matrix}$&$\begin{matrix}\text{-}2.0\text{e--}2\\2.2\text{e--}2\end{matrix}$&$\begin{matrix}\text{-}2.0\text{e--}2\\2.2\text{e--}2\end{matrix}$\\ \cline{2-7}
		\multicolumn{1}{|c|}{\multirow{5}{*}{}} & $6.4\text{e--}3$&$\begin{matrix}7.6\text{e--}3\\8.8\text{e--}3\end{matrix}$&$\begin{matrix}\text{-}4.3\text{e--}3\\8.4\text{e--}3\end{matrix}$&$\begin{matrix}\text{-}5.0\text{e--}3\\6.0\text{e--}3\end{matrix}$&$\begin{matrix}\text{-}5.0\text{e--}3\\5.7\text{e--}3\end{matrix}$&$\begin{matrix}\text{-}5.0\text{e--}3\\5.7\text{e--}3\end{matrix}$\\ \cline{2-7}
		\multicolumn{1}{|c|}{\multirow{5}{*}{}} & $1.6\text{e--}3$&$\begin{matrix}1.7\text{e--}2\\2.5\text{e--}3\end{matrix}$&$\begin{matrix}3.3\text{e--}3\\4.3\text{e--}3\end{matrix}$&$\begin{matrix}\text{-}1.2\text{e--}3\\2.3\text{e--}3\end{matrix}$&$\begin{matrix}\text{-}1.3\text{e--}3\\1.5\text{e--}3\end{matrix}$&$\begin{matrix}\text{-}1.3\text{e--}3\\1.4\text{e--}3\end{matrix}$\\ \cline{2-7}
		\multicolumn{1}{|c|}{\multirow{5}{*}{}} & $4.0\text{e--}4$&$\begin{matrix}2.0\text{e--}2\\6.5\text{e--}4\end{matrix}$&$\begin{matrix}1.2\text{e--}2\\1.7\text{e--}3\end{matrix}$&$\begin{matrix}3.7\text{e--}4\\1.6\text{e--}3\end{matrix}$&$\begin{matrix}\text{-}3.1\text{e--}4\\5.9\text{e--}4\end{matrix}$&$\begin{matrix}\text{-}3.1\text{e--}4\\3.7\text{e--}4\end{matrix}$\\ \cline{2-7}
		\multicolumn{1}{|c|}{\multirow{5}{*}{}} & $0\text{e+}0$&$\begin{matrix}2.6\text{e--}1\\0\text{e+}0\end{matrix}$&$\begin{matrix}2.6\text{e--}1\\0\text{e+}0\end{matrix}$&$\begin{matrix}2.6\text{e--}1\\0\text{e+}0\end{matrix}$&$\begin{matrix}2.6\text{e--}1\\0\text{e+}0\end{matrix}$&$\begin{matrix}2.6\text{e--}1\\0\text{e+}0\end{matrix}$\\ \hline
	\end{tabular}
\end{table}

\abbTab~\ref{tab:Exp_OCP_iter} shows the sum of the number of all inner iterations of MALM and QPM for respective $N,\pval$. We see the same trend as for the circle problem: QPM converges in a few iterations when $\pval$ is moderate. In contrast, when $\varepsilon,\pval$ both decrease then its iteration count increases. In contrast, MALM converges reliably for all $N,\pval$ in the upper right triangle, including those where $N$ is very large and $\pval$ very small. The last row shows that ALM does not convergence (n.c.) within $150$ iterations for any mesh size.

\begin{table}[tb]
\caption{Number of iterations for MALM/ALM and QPM for the Optimal Control Problem w.r.t.\ $N,\pval$. Fewer iterations mean better computational efficiency. MALM with $\omega=0$ (i.e.\ ALM) does not converge (n.c.) for this problem. QPM is not applicable (n.a.) when $\omega=0$.}
	\label{tab:Exp_OCP_iter}
	\centering
	\begin{tabular}{cc||c|c|c|c|c|}    \cline{3-7}\multicolumn{1}{l}{}&&\multicolumn{5}{c|}{$N$}\\ \cline{2-7}  
		\multicolumn{1}{l|}{}&$\begin{matrix}\#_\text{MALM}\\\#_\text{QPM}\end{matrix}$& $16$& $64$& $256$& $1024$& $4096$\\ \hline\hline
		\multicolumn{1}{|c|}{\multirow{5}{*}{$\pval$}} & $1.0\text{e--}1$& $\begin{matrix}\text{16}\\ \text{7}\end{matrix}$& $\begin{matrix}\text{16}\\ \text{7}\end{matrix}$& $\begin{matrix}\text{14}\\ \text{7}\end{matrix}$& $\begin{matrix}\text{12}\\ \text{7}\end{matrix}$& $\begin{matrix}\text{10}\\ \text{7}\end{matrix}$\\ \cline{2-7}
		\multicolumn{1}{|c|}{\multirow{5}{*}{}} & $2.5\text{e--}2$& $\begin{matrix}\text{22}\\ \text{7}\end{matrix}$& $\begin{matrix}\text{19}\\ \text{12}\end{matrix}$& $\begin{matrix}\text{16}\\ \text{10}\end{matrix}$& $\begin{matrix}\text{14}\\ \text{10}\end{matrix}$& $\begin{matrix}\text{13}\\ \text{10}\end{matrix}$\\ \cline{2-7}
		\multicolumn{1}{|c|}{\multirow{5}{*}{}} & $6.4\text{e--}3$& $\begin{matrix}\text{23}\\ \text{11}\end{matrix}$& $\begin{matrix}\text{20}\\ \text{13}\end{matrix}$& $\begin{matrix}\text{19}\\ \text{20}\end{matrix}$& $\begin{matrix}\text{18}\\ \text{19}\end{matrix}$& $\begin{matrix}\text{14}\\ \text{19}\end{matrix}$\\ \cline{2-7}
		\multicolumn{1}{|c|}{\multirow{5}{*}{}} & $1.6\text{e--}3$& $\begin{matrix}\text{25}\\ \text{12}\end{matrix}$& $\begin{matrix}\text{23}\\ \text{9}\end{matrix}$& $\begin{matrix}\text{20}\\ \text{42}\end{matrix}$& $\begin{matrix}\text{20}\\ \text{67}\end{matrix}$& $\begin{matrix}\text{17}\\ \text{101}\end{matrix}$\\ \cline{2-7}
		\multicolumn{1}{|c|}{\multirow{5}{*}{}} & $4.0\text{e--}4$& $\begin{matrix}\text{28}\\ \text{16}\end{matrix}$& $\begin{matrix}\text{27}\\ \text{19}\end{matrix}$& $\begin{matrix}\text{23}\\ \text{30}\end{matrix}$& $\begin{matrix}\text{23}\\ \text{68}\end{matrix}$& $\begin{matrix}\text{20}\\ \text{60}\end{matrix}$\\ \cline{2-7}
		\multicolumn{1}{|c|}{\multirow{5}{*}{}} & $0\text{e+}0$& $\begin{matrix}\text{n.~c.}\\ \text{n.~a.}\end{matrix}$& $\begin{matrix}\text{n.~c.}\\ \text{n.~a.}\end{matrix}$& $\begin{matrix}\text{n.~c.}\\ \text{n.~a.}\end{matrix}$& $\begin{matrix}\text{n.~c.}\\ \text{n.~a.}\end{matrix}$& $\begin{matrix}\text{n.~c.}\\ \text{n.~a.}\end{matrix}$\\ \hline
	\end{tabular}
\end{table}

\subsection{Interpretation of the Results}
For this test problem, \abbtab~\ref{tab:Exp_OCP_conv} demonstrates that the numerical solution converges to the optimal control solution when $h,\pval\searrow+0$ together.

\abbTab~\ref{tab:Exp_OCP_iter} shows that QPM converges fast when $N$ is small and $\pval$ is moderate. In contrast to this, for large $N$ and small $\pval$, MALM is clearly more efficient. However, \abbtab~\ref{tab:Exp_OCP_conv} reveals that large $N$ and small $\pval$ are a necessity for the numerical computation of accurate optimal control solutions.

Importantly, both experiments make clear that ALM is unsuitable for solving applications with inconsistent constraints: For the experiment depicted in Fig.~\ref{fig:example1geometry}, ALM converges to $\bx_A$ (white star) whenever $\varepsilon>0$. For optimal control problems, the magnitude of $\varepsilon$ models consistency errors of discretizations, in which case $\bx_B$ is the sought solution of a well-conditioned problem, whereas $\bx_A$ is the unsought solution of an ill-conditioned problem. \abbtab~\ref{tab:Exp_OCP_iter} shows that ALM fails to converge for the discretized control problem because it attempts to seek $\bx_A$, which is numerically hard. \abbtab~\ref{tab:Exp_OCP_conv} shows further that the exact minimizer is undesired here because $\pval=0$ results in a bad balance between the goals of minimizing both $r$ and $\delta J$, illustrated in Fig.~\ref{fig:ocp_study}.

\section{Conclusions}
\balance
We derived a modified version of ALM, called MALM. MALM outperforms QPM when minimizing unconstrained quadratic penalty programs \eqref{eqn:UCQPP} when $\pval$ is very small, in a similar manner as  ALM outperforms QPM when solving equality-constrained programs~\eqref{eqn:ECP}. 

The efficiency of MALM for the minimization of quadratic penalty functions has been demonstrated with numerical experiments. These experiments show that there are problem instances where it is beneficial to solve a problem of type \eqref{eqn:UCQPP} rather than \eqref{eqn:ECP}, one important class arising from integral penalty discretizations of optimal control problems.

In this paper we have presented the method in isolated form for ``approximately'' equality constrained programs in the sense that $c(\bx)\approx\bO$. Future work could extend the approach to problems with both equality and inequality constraints. Extensions of ALM for inequality constraints have been proposed in \cite{Powell2,Rockafellar1}, which can form a basis for similar extensions of MALM.



\bibliographystyle{plain} 		
\bibliography{MALMconf_refs_Martin_Eric}

\end{document}